\documentclass[12pt,a4paper]{amsart}
\usepackage{amsfonts,color}
\usepackage{amsthm}
\usepackage{amsmath}
\usepackage{amscd}
\usepackage{amssymb}
\usepackage[latin2]{inputenc}
\usepackage{t1enc}
\usepackage[mathscr]{eucal}
\usepackage{indentfirst}
\usepackage{graphicx}
\usepackage{graphics}
\usepackage{pict2e}
\usepackage{epic}
\usepackage{url}
\usepackage{epstopdf}
\usepackage{comment}
\usepackage{tikz}
\usepackage{xcolor}
\usepackage{fancyhdr}

\numberwithin{equation}{section}
\usepackage[margin=2.6cm]{geometry}

\allowdisplaybreaks

\theoremstyle{plain}
\newtheorem{Th}{Theorem}[section]

 \theoremstyle{definition}

\newtheorem{Rem}[Th]{Remark}
\newtheorem{?}[Th]{Problem}

\begin{document}

\title{Short proof of a theorem of Brylawski on the coefficients of the Tutte polynomial}

\author{Csongor Beke }

\address{Trinity College Cambridge, CB2 1TQ, United Kingdom \\
}

\email{bekecsongor@gmail.com}

\author{Gergely K\'al Cs\'aji}

\address{ELTE: E\"{o}tv\"{o}s Lor\'{a}nd University \\ H-1117 Budapest
\\ P\'{a}zm\'{a}ny P\'{e}ter s\'{e}t\'{a}ny 1/C  \and Institute of Economics, Centre for Economic and Regional Studies, Hungarian Academy of
Sciences, H-1097 Budapest, T\'{o}th K\'{a}lm\'{a}n u. 4\\ 
}

\email{csaji.gergely@gmail.com}

\author[P. Csikv\'ari]{P\'{e}ter Csikv\'{a}ri}

\address{Alfr\'ed R\'enyi Institute of Mathematics, H-1053 Budapest Re\'altanoda utca 13/15 \and ELTE: E\"{o}tv\"{o}s Lor\'{a}nd University \\ Mathematics Institute, Department of Computer
Science \\ H-1117 Budapest
\\ P\'{a}zm\'{a}ny P\'{e}ter s\'{e}t\'{a}ny 1/C}

\email{peter.csikvari@gmail.com}

\author{S\'ara Pituk}

\address{ELTE: E\"{o}tv\"{o}s Lor\'{a}nd University \\ H-1117 Budapest
\\ P\'{a}zm\'{a}ny P\'{e}ter s\'{e}t\'{a}ny 1/C
\\}

\email{pituksari@gmail.com}

\thanks{P\'eter Csikv\'ari is supported by the MTA-R\'enyi Counting in  Sparse Graphs ``Momentum'' Research Group. Csongor Beke and S\'ara Pituk are partially supported by the same research group. Gergely K\'al Cs\'aji is partially supported by the MTA-ELTE Matroid Optimization ``Momentum'' Research Group.}



\begin{abstract} 
In this short note we show that a system $M=(E,r)$ with a ground set $E$ of size $m$ and (rank) function $r: 2^E\to \mathbb{Z}_{\geq 0}$ satisfying $r(S)\leq \min(r(E),|S|)$ for every set $S\subseteq E$, the Tutte polynomial 
$$T_M(x,y):=\sum_{S\subseteq E}(x-1)^{r(E)-r(S)}(y-1)^{|S|-r(S)},$$
written as $T_M(x,y)=\sum_{i,j}t_{ij}x^iy^j$, satisfies that for any integer $h \geq 0$, we have
$$\sum_{i=0}^h\sum_{j=0}^{h-i}\binom{h-i}{j}(-1)^jt_{ij}=(-1)^{m-r}\binom{h-r}{h-m},$$
where $r=r(E)$, and we use the convention that when $h<m$, the binomial coefficient $\binom{h-r}{h-m}$ is interpreted as $0$. 

This generalizes a theorem of Brylawski on matroid rank functions and $h<m$, and a theorem of Gordon for $h\leq m$ with the same  assumptions on the rank function.

The proof presented here is significantly shorter than the previous ones. We only use the fact that the Tutte polynomial $T_M(x,y)$ simplifies to $(x-1)^{r(E)}y^{|E|}$ along the hyperbola $(x-1)(y-1)=1$.

\end{abstract}

\maketitle

\section{Introduction}

For a graph $G=(V,E)$ with $v(G)$ vertices and $e(G)$ edges, the Tutte polynomial $T_G(x,y)$ is defined as
$$T_G(x,y)=\sum_{A\subseteq E}(x-1)^{k(A)-k(E)}(y-1)^{k(A)+|A|-v(G)},$$
where $k(A)$ denotes the number of connected components of the graph $(V,A)$, see \cite{tutte1954contribution}. There are many excellent surveys about the properties of the Tutte polynomial and its applications \cite{brylawski1992tutte,crapo1969tutte,ellis2011graph,welsh1999tutte}. 

In this paper, we concentrate on Brylawski's identities concerning the Tutte polynomial. Written as a usual bivariate polynomial $T_G(x,y)=\sum_{i,j}t_{ij}x^iy^j$, the coefficients $t_{ij}$
encode the number of certain spanning trees, namely spanning trees with internal activity $i$ and external activity $j$ with respect to a fixed ordering of the edges, for details see \cite{tutte1954contribution}. It is not hard to prove that $t_{00}=0$ and $t_{10}=t_{01}$ if the graph $G$ has at least $2$ edges. In general, Brylawski \cite{brylawski1992tutte} proved that a collection of linear relations hold true between the coefficients of the Tutte polynomial.
Namely, he proved that if $0\leq h<e(G)$, then 
$$\sum_{i=0}^h\sum_{j=0}^{h-i}\binom{h-i}{j}(-1)^jt_{ij}=0.$$
In particular, the third relation gives that if $e(G)\geq 3$, then $t_{20}-t_{11}+t_{02}=t_{10}.$
Note that Brylawski \cite{brylawski1992tutte} proved these identities not only for the Tutte polynomial of a graph, but for the Tutte polynomial of an arbitrary matroid $M$. The Tutte polynomial $T_M(x,y)$ of a matroid $M=(E,r)$ is defined by
$$T_M(x,y)=\sum_{S\subseteq E}(x-1)^{r(E)-r(S)}(y-1)^{|S|-r(S)},$$
where $r(S)$ is the rank of a set $S\subseteq E$. The Tutte polynomial of a graph $G$ simply corresponds to the cycle matroid $M$ of the graph $G$. Note that the rank function $r:2^E\to \mathbb{Z}_{\geq 0}$ of a matroid satisfies the following axioms:\\
(R1) for any $A\subseteq E$ we have $r(A)\leq |A|$,\\
(R2) (submodularity) for any $A,B\subseteq E$ we have
$$r(A\cap B)+r(A\cup B)\leq r(A)+r(B),$$
(R3) (monotonicity) for any $A\subseteq E$ and $x\in E$ we have
$$r(A)\leq r(A\cup \{x\})\leq r(A)+1.$$
Gordon \cite{gordon2014tutte} calls a function $r:2^E\to \mathbb{Z}_{\geq 0}$ a rank function on a ground set $E$ if it satisfies $r(A)\leq \min(r(E),|A|)$ for every set $A\subseteq E$. He showed that for a system $M=(E,r)$ the coefficients of $T_M(x,y)$ satisfy Brylawski's identities if $r$ is a rank function without the assumptions of submodularity and monotonicity. He also extended Brylawski's identities to the case $h=|E|$.

Here we extend the work of Gordon and Brylawski for $h>|E|$, and also simplify the proof significantly. We only use the special form of the polynomial, namely that it simplifies to $(x-1)^{r(E)}y^{|E|}$  along the hyperbola $(x-1)(y-1)=1$. We use  exactly the same assumptions on the function $r:2^E\to \mathbb{Z}_{\geq 0}$ as Gordon.  Our generalized Brylawski's identities are the following.

\begin{Th}[Generalized Brylawski's identities]\label{matr} Let $M=(E,r)$, where $E$ is a set, and  $r: 2^E\to \mathbb{Z}_{\geq 0}$ is a  function on the subsets of $E$ satisfying $r(S)\leq \min(r(E),|S|)$ for every set $S\subseteq E$.
Let
$$T_M(x,y)=\sum_{S\subseteq E}(x-1)^{r(E)-r(S)}(y-1)^{|S|-r(S)}$$
be the Tutte polynomial of the system $M=(E,r)$.
Let $m$ denote the size of $E$, and let $r=r(E)$. The coefficients $t_{ij}$ of Tutte polynomial  $T_M(x,y)=\sum_{i,j}t_{ij}x^iy^j$ satisfy the following identities. For any integer $h \geq 0$, we have
$$\sum_{i=0}^h\sum_{j=0}^{h-i}\binom{h-i}{j}(-1)^jt_{ij}=(-1)^{m-r}\binom{h-r}{h-m},$$
with the convention that when $h<m$, the binomial coefficient $\binom{h-r}{h-m}$ is interpreted as $0$.
\end{Th}

In particular, by specializing Theorem~\ref{matr} for the cycle matroid of a graph $G$ we get the following.

\begin{Th}[Generalized Brylawski's identities for graphs]\label{graph}
Let $G$ be any graph with $n$ vertices, $m$ edges and $c$ connected components. Let $T_G(x,y)=\sum_{i,j}t_{ij}x^iy^j$ be the Tutte polynomial of the graph $G$. Then for any integer $h\geq 0$, we have
$$\sum_{i=0}^h\sum_{j=0}^{h-i}\binom{h-i}{j}(-1)^jt_{ij}=(-1)^{m-n+c}\binom{h-n+c}{h-m},$$
with the convention that when $h<m$, the binomial coefficient $\binom{h-n+c}{h-m}$ is interpreted as $0$.
\end{Th}
\bigskip

%\noindent \textbf{The paper is organized as follows.} In the next section, we prove Theorem~\ref{matr}.

\section{Proof of Theorem~\ref{matr}}
This entire section is devoted to the proof of Theorem~\ref{matr}.

Let $r=r(E)$ and $m=|E|$. By definition,
$$T_M(x,y)=\sum_{S\subseteq E}(x-1)^{r(E)-r(S)}(y-1)^{|S|-r(S)}.$$
Let us introduce a new variable $z$, and plug in $x=\frac{z}{z-1}$ and $y=z$. Then
$$T_M\left(\frac{z}{z-1}, z\right)=\sum_{S\subseteq E}(z-1)^{|S|-r}=(z-1)^{-r}z^{m}=\frac{z^{m}}{(z-1)^{r}}.$$
Since $T_M(x,y)=\sum_{i,j}t_{ij}x^iy^j$, we have
$$T_M\left(\frac{z}{z-1}, z\right)=\sum_{i,j}t_{ij}\left(\frac{z}{z-1}\right)^i z^j=\frac{z^m}{(z-1)^{r}}.$$
Hence
$$\sum_{i,j}t_{i,j}z^{i+j}(z-1)^{r-i}=z^m.$$
Note that if $i> r$, then $t_{ij}=0$ as $r(S)\geq 0$ for every set $S$.  
Hence, both sides are polynomials in $z$, so we can compare the coefficients of $z^k$. 
\begin{align}
\sum_{i,j}t_{i,j}(-1)^{r-k+j}\binom{r-i}{k-(i+j)}=\delta_{k,m},
\end{align}
where $\delta_{k,m}$ is $1$ if $k=m$, and $0$ otherwise.
This is not yet exactly Brylawski's identity, but taking appropriate linear combinations of these equations yields Brylawski's identities. Let
$$C_{h,k}=(-1)^k\binom{h-r}{h-k}.$$
Then
$$\sum_{k=0}^hC_{h,k}\left(\sum_{i,j}t_{i,j}(-1)^{r-k+j}\binom{r-i}{k-(i+j)}\right)=C_{h,m}.$$
Then
\begin{align*}
C_{h,m}&=\sum_{k=0}^hC_{h,k}\left(\sum_{i,j}t_{i,j}(-1)^{r-k+j}\binom{r-i}{k-(i+j)}\right) \\
&=\sum_{k=0}^h(-1)^k\binom{h-r}{h-k}\left(\sum_{i,j}t_{i,j}(-1)^{r-k+j}\binom{r-i}{k-(i+j)}\right)\\
&=\sum_{i,j}t_{i,j}(-1)^{r+j}\left(\sum_{k=0}^h\binom{h-r}{h-k}\binom{r-i}{k-(i+j)}\right)\\
&=\sum_{i,j}t_{i,j}(-1)^{r+j}\binom{h-i}{h-(i+j)}\\
&=\sum_{i,j}\binom{h-i}{j}t_{i,j}(-1)^{r+j}.
\end{align*}
Hence
$$\sum_{i,j}\binom{h-i}{j}t_{i,j}(-1)^{j}=(-1)^{m-r}\binom{h-r}{h-m}.$$

\begin{Rem} Once one conjectures Theorem~\ref{graph}, then it can be proved by the deletion-contraction identities via simple induction on $h$ even for matroids. The more general Theorem~\ref{matr} can be proved by certain recursions akin to deletion-contraction too, as  was shown by Gordon \cite{gordon2014tutte}, but seems to be considerably more work than the proof presented in this paper.
\end{Rem}

\noindent \textbf{Acknowledgment.} We are very grateful to the anonymous reviewers for the careful reading and the suggestions leading to a significant improvement in the presentation of this paper.
\bigskip

\noindent \textbf{Conflict of interest.} All authors declare that they have no conflicts of interest.

\bibliography{references}
\bibliographystyle{plain}

\end{document}